\documentclass[12pt]{article}

\usepackage{palatino,url}
\begin{document}

\section{Title}
AN ANALYTIC PROOF OF FERMAT-WILES THEOREM
\section{author}
Jamel Ghanouchi 

 ghanouchi.jamel@planet.tn
\section{Abstract}
 11 pages (MSC=11)

The present proof begins by changing the data of the problem. In
stead of an only equation, the traditional :
$$U^n=X^n+Y^n$$
We define two equations :
$$u^n=x^n+y^n$$
and
$$\frac{1}{z^n}=\frac{1}{x^n}+\frac{1}{y^n}$$
which allows to define four decreasing and restricted sequences.
We study those sequences and an elementary calculus allows us to
conclude quickly that there is no solution for $n>2$.

\section{Preliminary}

Our goal is to prove that if U, X, Y, n are positive integers,
they never verify for $n>2$ an equation such as :
(a)
$$U^n=X^n+Y^n$$
Our proof is supported by a reasoning that must conduct to an
absurd result. We will suppose that $\exists{U, X, Y}$ and n
positive integers which verify : (a).
If we choose u, x, y, z and
n positive integers which are equal to :
\\
$$u=U^2$$
\\
$$x=UX$$
\\
$$y=UY$$
\\
$$z=XY$$
\\
Then they verify :
\\
(1) $$u^n=(U^2)^n=U^nU^n=U^n(X^n+Y^n)=U^nX^n+U^nY^n=x^n+y^n$$
\\
and :
\\
(2)
$$\frac{1}{z^n}=\frac{1}{X^nY^n}=\frac{U^{2n}}{(UX)^n(UY)^n}$$
$$=\frac{u^n}{x^ny^n}=\frac{x^n+y^n}{x^ny^n}$$
$$=\frac{1}{x^n}+\frac{1}{y^n}$$
\\
\section{The proof}
Then if U, X, Y, n are positive integers which verify (a),
u, x, y, z, n as defined must verify simultaneously (1) and (2) :
\\
(1) $$u^n=x^n+y^n$$
\\
and (2)
\\
$$\frac{1}{z^n}=\frac{1}{x^n}+\frac{1}{y^n}$$ \\
\\
 Then :
\\
$$z^2=X^2Y^2$$
\\
$$u^n=(X^n+Y^n)^2$$
\\
$$x^n=X^n(X^n+Y^n)$$
\\
$$y^n=Y^n(X^n+Y^n)$$
\\
$p=GCD(X,Y)=HCF(X,Y)$
\\
$$X^2=p^2v_2^2$$
\\
and
\\
$$Y^2=p^2t_2^2$$
\\
with
\\
$$GCD(v_2,t_2)=HCF(v_2,t_2)=1$$
\\
$$z^2=p^4v_2^2t_2^2$$
\\
$$z=p^2v_2t_2$$
\\
$$x^n=p^{2n}(v_2^n)(v_2^n+t_2^n)$$
\\
$$y^n=p^{2n}(t_2^n)(v_2^n+t_2^n)$$
\\
$$w_2^n=v_2^n+t_2^n$$
\\
with
\\
$$u^n=p^{2n}w^n=p^{2n}w_2^{2n}$$
\\
For
$n=1$
\\
$$w_2=v_2+t_2$$
\\
has an infinity of solutions
\\
for $n=2$
\\
$$w=w_2^2=v_2^2+t_2^2=w_2^2$$
\\
has also an infinity of solutions and we continue the reasoning
for $n>2$
\\
$v_2^n+t_2^n=w_2^n$
 $v_2^n$
 and
 $t_2^n$
 are coprime two by two.
\\
$\forall{v_2,t_2}$ positive integers, $\exists{z_2}$ positive
real number which verifies
\\
$$\frac{1}{z_2^n}=\frac{1}{v_2^n}+\frac{1}{t_2^n}$$
\\
and :
\\
$$(v_2^n+t_2^n)z_2^n=v_2^nt_2^n$$
\\
which means
\\
$$v_2^n(t_2^n-z_2^n)=z_2^nt_2^n$$
\\
We pose :
\\
$$t_3^n=t_2^n-z_2^n=\frac{z_2^nt_2^n}{v_2^n}$$
\\
and :
$$t_2^n(v_2^n-z_2^n)=z_2^nv_2^n$$
\\
We pose :
\\
$$v_3^n=v_2^n-z_2^n=\frac{z_2^nv_2^n}{t_2^n}$$
\\
and :
$$v_3^nt_3^n=z_2^{2n}$$
\\
Then:
\\
$$v_2^n=v_3^n+z_2^n=v_3^n+\sqrt{v_3^nt_3^n}$$
\\
and \\
$$t_2^n=t_3^n+z_2^n=t_3^n+\sqrt{v_3^nt_3^n}$$
\\
and \\
$$v_2^n+t_2^n=(\sqrt{v_3^n}+\sqrt{t_3^n})^2>v_3^n+t_3^n>1$$
\\
and $v_2^n+t_2^n$ integer
\\
$$v_2^n=\sqrt{v_3^n}(\sqrt{v_3^n}+\sqrt{t_3^n})>v_3^n>0$$
\\
and $v_2^n$ integer
\\
$$t_2^n=\sqrt{t_3^n}(\sqrt{v_3^n}+\sqrt{t_3^n})>t_3^n>0$$
\\
and $t_2^n$ integer
\\
$$z_2^n=\frac{v_2^nt_2^n}{v_2^n+t_2^n}=\sqrt{v_3^nt_3^n}>z_3^n=\frac{v_3^nt_3^n}{w_3^n}>0$$
\\
and $z_2^n$  real
\\
and :
\\
$\forall{v_3,t_3}$
positive real numbers,
$\exists{z_3}$
positive
real number which verifies
\\
$$\frac{1}{z_3^n}=(\frac{1}{v_3^n})+(\frac{1}{t_3^n})$$
\\
$v_3$
 $t_3$
  $z_3$
  are real numbers, until infinity ! For i :
\\
$$v_i^n+t_i^n=(\sqrt{v_{i+1}^n}+\sqrt{t_{i+1}^n})^2>v_{i+1}^n+t_{i+1}^n>1$$
\\
and :
$v_i^n+t_i^n$
is a real number for
$i>2$
\\
$$v_i^n=\sqrt{v_{i+1}^n}(\sqrt{v_{i+1}^n}+\sqrt{t_{i+1}^n})>v_{i+1}^n>0$$
\\
$v_i^n$
is a real number for
$i>2$
\\
$$t_i^n=\sqrt{t_{i+1}^n}(\sqrt{v_{i+1}^n}+\sqrt{t_{i+1}^n})>t_{i+1}^n>0$$
\\
$t_i^n$
is a real number for
$i>2$
\\
$$z_i^n=\frac{v_i^nt_i^n}{v_i^n+t_i^n}=\sqrt{v_{i+1}^nt_{i+1}^n}>z_{i+1}^n=\frac{(v_{i+1}^n)(t_{i+1}^n)}{v_{i+1}^n+t_{i+1}^n}>0$$
\\
$z_i^n$
is a real number for
$i>2$
and :
\\
$$\frac{1}{z_{i+1}^n}=\frac{1}{v_{i+1}^n}+\frac{1}{t_{i+1}^n}$$
\\
With this recurrence :
\\
$$v_i^n+t_i^n=(\sqrt{v_{i+1}^n}+\sqrt{t_{i+1}^n})^2$$
\\
Then :
\\
$$v_i^n=v_2^{n{2^{i-2}}}\prod_{j=0}^{j={i-3}}{(v_2^{n{2^j}}+t_2^{n{2^j}})^{-1}}$$
\\
which is (H)
\\
$$t_i^n=t_2^{n{2^{i-2}}}\prod_{j=0}^{j={i-3}}{(v_2^{n{2^j}}+t_2^{n{2^j}})^{-1}}$$
\\
which is (H'), because
\\
$$v_3^n=\frac{v_2^{2n}}{v_2^n+t_2^n}$$
\\
it is verified ! We suppose that we have : (H)and (H')
\\
$$v_i^n=\sqrt{v_{i+1}^n}(\sqrt{v_{i+1}^n}+\sqrt{t_{i+1}^n})=\sqrt{v_{i+1}^n}(v_i^n+t_i^n)^{\frac{1}{2}}$$
\\
Then :
\\
$$v_{i+1}^n=v_i^{2n}{(v_i^n+t_i^n)}^{-1}$$
\\
But (H) and (H')
\\
$$v_{i+1}^n={v_2^{2^{i-1}n}}\prod_{j=0}^{j=i-2}{(v_2^{{2^j}n}+t_2^{{2^j}n})^{-1}}$$
\\
It is proved ! Now :
\\
$$\forall{i>2}$$
\\
$$v_i^n-t_i^n=(v_2^n-\sum_{j=3}^{j=i}{\sqrt{v_j^nt_j^n}})-(t_2^n-\sum_{j=3}^{j=i}{\sqrt{v_j^nt_j^n}})=v_2^n-t_2^n$$
\\
but $\forall{a,b}$ \\
$$a-b=(a^{2^{i-2}}-b^{2^{i-2}})\prod_{j=0}^{j={i-3}}{(a^{2^j}+b^{2^j})^{-1}}$$
\\
which is (E). Then \\
$$v_2^n-t_2^n=(v_2^{n{2^{i-2}}}-t_2^{n{2^{i-2}}})\prod_{j=0}^{j={i-3}}{(v_2^{n{2^j}}+t_2^{n{2^j}})^{-1}}$$
\\
which is (E') We can also write with
\\
$$v_2\neq{t_2}$$
\\
then
\\
$$v_i\neq{t_i}$$
\\
$$v_i^n=\frac{v_2^{n{2^{i-2}}}}{v_2^{n{2^{i-2}}}-t_2^{n{2^{i-2}}}}(v_2^n-t_2^n)$$
\\
and
\\
$$t_i^n=\frac{t_2^{n{2^{i-2}}}}{v_2^{n{2^{i-2}}}-t_2^{n{2^{i-2}}}}(v_2^n-t_2^n)$$
\\
Also
\\
$$v_i^n=\frac{v_i^n}{v_i^n-t_i^n}v_2^n(1-\frac{t_2^n}{v_2^n})=\frac{v_i^n}{v_i^n-t_i^n}v_2^n(1-\frac{t_i^{\frac{n}{2^{i-2}}}}{v_i^{\frac{n}{2^{i-2}}}})$$
\\
then
\\
$$v_2^n=\frac{v_i^{\frac{n}{2^{i-2}}}}{v_i^{\frac{n}{2^{i-2}}}-t_i^{\frac{n}{2^{i-2}}}}(v_2^n-t_2^n)$$
\\
and also \\
$$t_2^n=\frac{t_i^{\frac{n}{2^{i-2}}}}{v_i^{\frac{n}{2^{i-2}}}-t_i^{\frac{n}{2^{i-2}}}}(v_2^n-t_2^n)$$
\\
$$w_2^n=v_2^n+t_2^n=\frac{v_i^{\frac{n}{2^{i-2}}}+t_i^{\frac{n}{2^{i-2}}}}{v_i^{\frac{n}{2^{i-2}}}-t_i^{\frac{n}{2^{i-2}}}}(v_2^n-t_2^n)$$
\\
then
\\
$$\forall{i>2}$$
\\
for $$n>2$$
\\
$$v_i^{\frac{n}{2^{i-2}}}-t_i^{\frac{n}{2^{i-2}}}=(v_2^n-t_2^n)\prod_{j=0}^{j={i-3}}{(v_i^{\frac{n{2^j}}{2^{i-2}}}+t_i^{\frac{n{2^j}}{2^{i-2}}})^{-1}}$$
\\
and from (H) and (H')
\\
$$v_i^{\frac{n}{2^{i-2}}}=v_2^n\prod_{j=0}^{j={i-3}}{(v_2^{n{2^j}}+t_2^{n{2^j}})^{-{\frac{1}{2^{i-2}}}}}$$
\\
and \\
$$t_i^{\frac{n}{2^{i-2}}}=t_2^n\prod_{j=0}^{j={i-3}}{(v_2^{n{2^j}}+t_2^{n{2^j}})^{-{\frac{1}{2^{i-2}}}}}$$
\\
$$w_2^n=v_2^n+t_2^n=(v_i^{\frac{n}{2^{i-2}}}+t_i^{\frac{n}{2{i-2}}})^n\prod_{j=0}^{j={i-3}}{(v_2^{n{2^j}}+t_2^{n{2^j}})^{\frac{1}{2^{i-2}}}}$$
$$=\frac{v_i^{\frac{n}{2^{i-2}}}+t_i^{\frac{n}{2^{i-2}}}}{v_i^{\frac{n}{2^{i-2}}}-t_i^{\frac{n}{2^{i-2}}}}(v_2^n-t_2^n)$$
$$=(v_i^{\frac{n}{2^{i-2}}}+t_i^{\frac{n}{2^{i-2}}})\prod_{j=0}^{j={i-3}}{(v_i^{\frac{n{2^j}}{2^{i-2}}}+t_i^{\frac{n{2^j}}{2^{i-2}}})}$$
\\
then
\\
$$\prod_{j=0}^{j={i-3}}{(v_i^{\frac{n{2^j}}{2^{i-2}}}+t_i^{\frac{n{2^j}}{2^{i-2}}})}=\prod_{j=0}^{j={i-3}}{(v_2^{n{2^j}}+t_2^{n{2^j}})^{\frac{1}{2^{i-2}}}}$$
\\
$$\forall{i>2}$$
\\
then \\
$$(\prod_{j=0}^{j={i-3}}{(v_2^{n{2^j}}+t_2^{n{2^j}})})^{\frac{2^{i-3}-1}{2^{i-3}}}=1$$
\\
then \\
$$v_2^n=0$$
\\
or \\
$$t_2^n=0$$
\\
it means that the only solutions are the trivial
\\
$$X=0$$
\\
or
\\
$$Y=0$$
\\
it is impossible !
\section{Conclusion}

The recurrence conducts us to an impossibility for $n>2$, the
initial hypothesis is not correct, we never have for $n>2$,
simultaneously (1) et (2),

which means that we never have (a).

Fermat-Wiles theorem is proved for U, X, Y positive integers.

But, if
 $\forall{ X, Y, U, n>2}$ positive integers

$$U^n\neq{X^n+Y^n}$$

then

$$(adf)^n\neq{(cbf)^n+(ebd)^n}$$

$\forall{a,b,c,d,e,f}$ positive integers, after dividing by
$(bdf)^n$,

$$\frac{a^n}{b^n}\neq{\frac{c^n}{d^n}+\frac{e^n}{f^n}}$$
$\forall{a,b,c,d,e,f}$ positive integers and $n>2$ !

The proof is general !

By :

Jamel GHANOUCHI

ghanouchi.jamel@planet.tn

\end{document}